\input amstex
\magnification=1200
\voffset-1in
\documentstyle{amsppt}
\NoBlackBoxes
\def\Z{\Bbb Z}

\def\Q{\Bbb Q}
\def\H{\Cal H}
\def\a{\alpha}
\def\l{\lambda}
\def\e{{\bold e}}
\def\m{{\bold m}}

\def\sltwo{{\frak s\frak l _2 }}
\def\slk{\frak{sl}_k}
\def\<{\langle}
\def\>{\rangle}
\def\Uk{U_q \frak s\frak l _k }
\def\U2{U_q \frak s\frak l _2 }
\def\o{\otimes}
\topmatter
\title Kazhdan-Lusztig polynomials and canonical basis
\endtitle
\date September  22, 1997 \enddate
\author Igor Frenkel, Mikhail Khovanov\\
 and  Alexander Kirillov, Jr. \endauthor
\address Dept. of Mathematics, Yale University, New Haven, CT
06520--8283, USA\endaddress
\address School of Math., Inst. Adv. Study, Olden Lane, Princeton, NJ,
08540, USA\endaddress 
\email mikhailk\@math.ias.edu\endemail
\address Math. Dept., MIT, Rm. 2--179,  Cambridge, MA 02139, USA
\endaddress
\email kirillov\@math.mit.edu\endemail
\leftheadtext{Frenkel, Khovanov and Kirillov}
\abstract
	In this paper we show that the Kazhdan--Lusztig polynomials
	(and, more generally, parabolic KL polynomials) for the group
	$S_n$ coincide with the coefficients of the canonical basis in
	$n$th tensor power of the fundamental representation of the
	quantum group $\Uk$. We also use known results about canonical
	bases for $\U2$ to get a new proof of recurrent
	formulas for KL polynomials for maximal parabolic subgroups
	(geometrically, this case corresponds to Grassmanians), due to
	Lascoux-Sch\"utzenberger and Zelevinsky.
\endabstract

\endtopmatter

\document

\head 1. Review of the theory of Kazhdan-Lusztig polynomials
\endhead

In this section, we review the theory of Kazhdan-Lusztig
polynomials. We will use their generalization to the parabolic case,
defined by Deodhar (see \cite{D}). For the sake of completeness and to
fix notations, we list the main definitions and results here,
referring the reader to the original papers for more details.
To avoid confusion with the theory of quantum groups, we will not use
the variable $q$ in the definition of the Hecke algebra; instead, we
will use $v=q^{-1/2}$. 

Let $W$ be a finite Weyl group with a set  of simple reflections
$S$. We denote by $l(w)$ the length of $w\in W$ with respect to the
generators $s\in S$. 

Let $\Cal H (W)$ be the Hecke algebra associated with $W$ (we will
usually omit $W$ and write just $\H$). By definition, it is an
associative algebra with unit over the field $\Q(v)$ with generators
$T_s, s\in S$ and relations 

$$\gathered
{\aligned T_s T_{s'}\dots =T_{s'}&T_s\dots\qquad\text{$n$ terms on each
	side }\\
	& \text{where $n$ is the order of $ss'$ in $W$}
\endaligned}\\
(T_s+1)(T_s-v^{-2})=0.
\endgathered
\tag 1.1
$$

For any subset $J\subset S$, let $W_J$ be the parabolic subgroup in
$W$ generated by $s\in J$. Denote by $W^J$ the set of minimal length
representatives of the cosets in $W/W_J$, and let $<$ be the
partial order on $W^J$ induced by the Bruhat order on $W$.

Let $M$ be a vector space over $\Q(v)$ with the basis $m_y, y\in
W^J$. 

\proclaim{Proposition 1.1\rm (see \cite{D, Lemmas~2.1, 2.2})}
Let $u$ be either $-1$ or $v^{-2}$. 

\roster\item The following formulas define the structure of an
$\H$-module on $M$:

$$T_sm_\sigma=\cases v^{-1}m_{s\sigma}+(v^{-2}-1)m_\sigma, \quad
			  & l(s\sigma)<l(\sigma),\\
		     v^{-1}m_{s\sigma}	
			  & l(s\sigma)>l(\sigma), \ s\sigma\in W^J\\
		     um_{\sigma}, 
			  &  l(s\sigma)>l(\sigma),\ s\sigma\notin W^J
		\endcases
\tag 1.2
$$	  
Note that if $\sigma\in W^J, l(s\sigma)< l(\sigma)$ then $s\sigma\in
W^J$. 
 			
\item If $\sigma \in W^J$ then $T_\sigma m_e=v^{-l(\sigma)}m_\sigma$.

\endroster
\endproclaim

\remark{Remarks}  1. Our notations slightly differ from those of
Deodhar: what we denote $m_\sigma$ in his notations would be
$q^{-l(\sigma)/2} m_\sigma$. Note also that there is a misprint in the
formula for $L(s)$ (immediately after the statement of Lemma~2.1) in
\cite{D}.

2. It is easy to see that $M$ is a deformation of the
induced representation $\text{Ind}_{W_J}^{W} 1$ of $W$. In particular,
if $J=\varnothing$ then $M$ is the left regular representation of $\H$
and the action is independent of  $u$. 
\endremark

Define an involution $\overline{\phantom{T}}:\Q(v)\to \Q(v)$ by
$f(v)\mapsto f(v^{-1})$. We will say that a map $\phi$ of
$\Q(v)$-vector spaces is antilinear if $\phi(f x)=\bar f \phi(x)$ for any
$f\in \Q(v)$ and any vector $x$.

We define antilinear involutions on $\H$ and $M$ by letting
$\overline{T_s}=T_s^{-1}=v^2T_s+v^2-1, \overline{m_e}=m_e$ and
$\overline{ab}=\bar a\bar b$. Note that one has 
$\overline{m_\sigma}=\overline{(v^{l(\sigma)}T_\sigma)}
m_e=m_\sigma+\sum_{\tau<\sigma}  c_{\tau\sigma}m_\tau$.

\proclaim{Theorem 1.2\rm (\cite{D, Proposition 3.2})}
Let us assume that we have fixed $u$ as in Proposition~\rom{1.1}. 
Then for every $\sigma\in W^J$ there exists a unique element $C_\sigma\in
M$ such that the following two conditions are satisfied:
$$
\aligned
\overline{C_\sigma}=&C_\sigma,\\
C_\sigma=&\sum\Sb \tau\in W^J\\\tau\le \sigma\endSb
		\a_{\tau\sigma}m_\tau,\\
&\quad \text{\rm where }\a_{\sigma\sigma}=1\quad 
	\text{\rm and } \a_{\tau\sigma}\in v^{-1}\Z[v^{-1}]
	\quad \text{\rm for }\tau<\sigma.
\endgathered 
\tag 1.3
$$

The elements $C_\sigma, \sigma\in W^J$ form a basis in $M$. 

\endproclaim

Note that the definition of $C_\sigma$ uses the involution 
$\overline{\phantom{T}}$, which was defined in terms of the action of
$\H$ and therefore depends on the choice of $u$. Thus, we have two
different antilinear involutions on the same space $M$ with fixed
basis $m_y$, which give rise to different bases $C_\sigma$.

Following Deodhar, we define {\bf parabolic Kazhdan-Lusztig
polynomials\/} $P^J_{\tau\sigma}\in \Z[q], \tau, \sigma\in W^J,
\tau\le\sigma$ by
$$
\a_{\tau\sigma}=(-v)^{l(\tau)-l(\sigma)} \overline{P^J_{\tau\sigma}}
\tag 1.4
$$
where $\a_{\tau\sigma}$ is defined by (1.3) and as before, $q=v^{-2}$.
This polynomials can be expressed via the usual KL polynomials as
follows:

\proclaim{Theorem 1.3 \rm{\cite{D, Proposition~3.4 and Remark~3.8}}}
\roster
\item For $u=-1$ we have
$$
P^J_{\tau\sigma}= P_{\tau w_J^0, \sigma w_J^0},
\tag 1.5
$$
where $w_J^0$ is the longest element in $W_J$ and $P_{y,w}$ is the
usual Kazhdan-Lusztig polynomial for $W$. 
\item For $u=v^{-2}$, 
$$
P^J_{\tau\sigma}=  \sum\Sb w\in W_J\\ \tau w\le\sigma\endSb
			(-1)^{l(w)} P_{\tau w, \sigma}.
\tag 1.6
$$
\endroster
\endproclaim

In this paper we will be interested in the case $W=S_n$. We
will consider $S_n$ as the group acting by permutations on the set of
all integer sequences  of length $n$, with the generators $s_i,
i=1,\dots, n-1$  acting in the standard way: $s_i(x_1, \dots,
x_n)=(x_1, \dots, x_{i+1}, x_i, \dots, x_n)$.

\head 2. Canonical basis for $\Uk$.
\endhead

In this section, we show that the parabolic Kazhdan-Lusztig basis for
$S_n$ can be obtained as a special case of the canonical basis for
representations of quantum groups. This is inspired by the results of 
\cite{GL}, where it is proved that the projectivization of the KL basis
coincides with the so-called ``special basis'' for quantum
groups. However, it would take us more time to define what a special
basis is and why it is the projectivization of the canonical
basis. For this reason, we are using the results of  \cite{GL} only as
motivation. 

Let us fix an integer $k\ge 2$. Let $\Uk$ be the quantum group
corresponding to the Lie algebra $\slk$, that is, an associative
algebra over the field $\Q(v)$ with generators $E_i, F_i, v^{H_i},
i=1,\dots , k-1$ and standard relations (see, for example, \cite{L},
where $v^{H_i}$ is denoted by $K_i$). We denote by $(\Uk)^{\pm}$
subalgebras in $\Uk$ generated by $E_i, v^{H_i}$ (respectively, by
$F_i, v^{H_i}$). We identify the weight lattice $P$ for $\slk$ with
$\Z^k/\Z\cdot(1,\dots, 1)$ by the rule $H_i(\l)=\l_i-\l_{i+1}, \l=(\l_1,
\dots, \l_k)\in
\Z^k$. 

Let $V$ be a $k$-dimensional vector space over $\Q(v)$ with the basis
$\e_1, \dots, \e_k$. We define on $V$ the structure of a representation
of $\Uk$ by 

$$\aligned
&E_i\e_{i+1}=\e_i, \quad E_i \e_j=0, \quad j\ne i+1,\\
&F_i\e_{i}=\e_{i+1}, \quad F_i \e_j=0, \quad j\ne i,\\
&v^{H_i} \e_j=\cases v \e_i, \quad & j=i,\\
		   v^{-1}\e_{i+1}, \quad & j=i+1,\\
		   \e_i, \quad & j\ne i, i+1
	    \endcases
\endaligned
\tag 2.1
$$

As is well-known (see, e.g., \cite{L}, \cite{CP}), $\Uk$ can be
endowed with the structure of a Hopf algebra with the
comultiplication

$$\aligned 
&\Delta E_i= E_i\o 1+ v^{H_i}\o E_i,\\
&\Delta F_i= F_i\o v^{-H_i}+1 \o F_i,\\
&\Delta v^{H_i}=v^{H_i}\o v^{H_i}
\endaligned 
\tag 2.2
$$ 
We do not
write formulas for the antipode and counit since we will not use them.

Moreover, $\Uk$ is quasitriangular: there is an element $\Cal R$
(universal $R$-matrix) in a certain completion of $(\Uk)^{\o 2}$,
which, among other properties, has the following one.  For any pair of
finite-dimensional representations $V, W$  with weight decomposition,
$R_{V\o W}=\Cal R|_{V\o W}$ is well-defined and $PR_{V\o W}:V\o W\to
W\o V$ is an isomorphism of $\Uk$-modules. In fact, this last property
uniquely defines $\Cal R$ if we also require that it have the
following form: 

$$
\Cal R=\bar C\bar \Theta
\tag 2.3$$
where 

$$\gathered
\bar C=v^{\sum x_a\o x_a},\\
\bar  \Theta = 1+\sum a_i\o b_i, \quad a_i\in (\Uk)^-, b_i\in (\Uk)^+
\endgathered\tag 2.4
$$
where $x_a$ is an orthonormal basis in the Cartan subalgebra of $\slk$
and $a_i$ (respectively, $b_i$) are of negative (respectively,
positive) weight. We refer the reader to \cite{L} or \cite{CP} for
details.  

\remark{Remark} The notations $\bar C, \bar \Theta$ are chosen to
agree with notations of Lusztig and in \cite{FK}.  Later we will define
the bar  involution $\overline{\phantom{T}}$ and elements $C, \Theta$ so
that $\bar C, \bar \Theta$ defined by (2.4) will indeed be the bar
conjugates of  $C, \Theta$.
\endremark

Let us consider the vector space $V^{\o n}$. This space has a basis
given by 

$$
\e_I=\e_{i_1}\o\dots\o \e_{i_n},\qquad
	I=(i_1,\dots, i_n)\in \{1,\dots, k\}^n.
\tag 2.5
$$

We define the action of $\Uk$ on $V^{\o n}$ using the comultiplication
$\Delta$. This space has a weight decomposition: if $\m\in \Z_+^k$ is
such that $\sum m_i=n$ then a basis in the weight subspace $V^{\o
n}[\m]$ is given by the vectors $\e_I$ such that every $a=1,\dots, k$
appears in the sequence $I=(i_1,\dots, i_n)$ exactly $m_a$ times. 

\proclaim{Proposition 2.1} 
\roster \item Let $P:V\o V\to V\o V$ be the permutation: $P(v\o w)=w\o
v$, and let $R_{V\o V}$ denote the action of the universal $R$-matrix
$\Cal R$ in $V\o V$.  Then

$$
PR_{V\o V}(\e_a\o \e_b)= v^{-1/n}
	\cases \e_b\o \e_a +(v-v^{-1}) \e_a\o \e_b,\quad& a<b,\\
		\e_b\o \e_a, \quad& a>b,\\
		v \e_a\o \e_b,\quad& a=b	
	\endcases
\tag 2.6
$$

\item 
Denote by $(PR)_i$ an endomorphism of $V^{\o n}$ which acts as
$PR_{V\o V}$ on the tensor product of $i$-th and $i+1$-th factors and
by identity on all other factors. Then the map 

$$
T_i\mapsto -v^{-1+1/n} (PR)_i
\tag 2.7
$$
defines an action of the Hecke algebra $\Cal H(S_n)$  on
$V^{\o n}$. This action commutes with the action of $\Uk$; in
particular, it preserves the weight subspaces.

\item 
For $\m\in \Z_+^k, \sum m_i=n$ let 

$$I^0(\m)=(\underbrace{k,\dots, k}_{\text{\rm $m_k$ times}}, \dots, 
	\underbrace{1,\dots, 1}_{\text{\rm $m_1$ times}}).
\tag 2.8
$$

Let $J\subset \{1,\dots, n-1\}$ be such that the corresponding
parabolic subgroup $W_J\subset S_n$ is the stabilizer of $I^0(\m)$,
i.e.  $W_J=S_{m_k}\times \dots \times S_{m_1}$, and let $M$ be the
corresponding vector space as defined in Section~\rom{1}. Then the map

$$\aligned
\phi_\m:M&\to V^{\o n}\\
	m_\sigma&\mapsto (-1)^{l(\sigma)}\e_{\sigma(I^0(\m))}
\endaligned
\tag 2.9$$
is an isomorphism of $\Cal H$ modules $M$ \rom{(}with $u=-1$\rom{)} and $V^{\o
n}[\m]$. 
\endroster
\endproclaim

\demo{Proof} Parts (1), (2) are due to Jimbo \cite{J} and by now are
well-known.  As for (3), it follows from the comparison of formulas
(1.2) and (2.6) and the following simple lemma, the proof of which is
left to the reader.

\proclaim{Lemma} Let $\m, J$ be as in the statement of the theorem,
and let $\sigma\in W^J$.  Let $a$ and $b$ be  $i$-th
and $i+1$-th entries of $\sigma(I^0(\m))$, respectively. Then
$l(s_i\sigma)>l(\sigma)$ if and only if $a\ge b$, and $s_i\sigma\in
W^J$ if and only if $a\ne b$. 
\endproclaim
\qed\enddemo

For future reference, we formulate the following result, which is very
close to the lemma above. 

\proclaim{Proposition 2.2} Let $\m, J$ be as in the previous
proposition. Then the map 
$$
\sigma\mapsto \sigma(I^0(\m))
$$
is a bijection between $W^J$ and the set of all sequences of weight
$\m$, and 
$$
l(\sigma)=\{I^0(\m)\}-\{\sigma(I^0(\m))\}
$$
where $\{I\}$ is the number of inversions in $I$, i.e. the number of
pairs of indices $a, b$ such that $a<b$ and $i_a>i_b$. 
\endproclaim

Recall that we define the bar involution on $\Q(v)$ by $\bar v=v^{-1}$
and we say that a map $\phi$ of $\Q(v)$ vector spaces is antilinear if
$\phi(f v)=\bar f\phi(v)$. Following Lusztig, we will also denote by
a bar the antilinear algebra automorphism of $\Uk$ defined on the
generators by

$$
\overline{E_i}=E_i,\quad 
 \overline{F_i}=F_i,\quad \overline{v^{H_i}}=v^{-H_i}.
$$

We extend the bar involution to $\Uk\o \Uk$ by the rule $\overline{x\o
y} =\bar x\o \bar y$ and denote by $C, \Theta$ the bar conjugates of
the elements $\bar C,\bar \Theta$ defined by (2.4). We will use the
following result due to Lusztig: $\Theta\bar\Theta=1$. Note also that
$C\bar C=1$ (obvious).

Let us define by induction an antilinear involution $\psi$ on tensor
powers of the module $V$ as follows:

\roster
\item  On $V$, $\psi$ is given by 

$$\psi\left(\sum a_i \e_i\right)=\sum \overline{a_i} \e_i.
\tag 2.10a$$

\item If $W_1, W_2$ are tensor powers of $V$ and the involution $\psi$
is already defined on $W_1, W_2$ then on $W_1\o W_2$ we have 

$$\psi(w_1\o w_2)=\Theta (\psi(w_1)\o \psi(w_2)). \tag 2.10b$$

As before, the action of $\Uk$ \rom{(}and thus, of $\Theta$\rom{)} on
$W_i$ is defined using the comultiplication $\Delta$. 

\endroster

It follows from the results of Lusztig that $\psi$ is well defined on
any tensor power of $V$ and that $\Psi^2=1$.

We will use the involution $\psi$ to define the canonical basis in
$V^{\o n}$. The following proposition is a special case of the
definition in \cite{L, Chapter 27}. 

\proclaim{Proposition 2.3} For every $I\in \{1, \dots, k\}^n$ there
exists a unique element $b_I\in V^{\o n}$ such that 

$$\gathered
\psi(b_I)=b_I,\\
b_I-\e_I\in \bigoplus_{I'} v^{-1}\Z[v^{-1}]\e_{I'},
\endgathered
\tag 2.11
$$
where the sum is taken over all $I'$ having the  same weight as $I$. 

The elements $b_I$ form a basis in $V^{\o n}$ which is called the
{\bf canonical basis}. 
\endproclaim

The first important  result of this section is the following proposition.

\proclaim{Proposition 2.4} The isomorphism $\phi_\m$ constructed in
Proposition~\rom{2.1} identifies the bar involution in $M$ \rom{(}for
$u=-1$\rom{)} with the involution $\psi$ in $V^{\o n}[\m]$.
\endproclaim
\demo{Proof} 
Since $M$ is spanned by $T_ym_e$, it suffices to prove that
$\psi(\phi_\m(m_e))=\phi_{\m}(m_e)$, $\psi
\phi_\m(T_i)=\phi_\m(T_i^{-1})\psi$, or, equivalently,
$\psi(\e_{I^0(\m)})=\e_{I^0(\m)}, \psi(PR)_i=(PR)_i^{-1}\psi$. 

To prove the first identity, let $<$ be the partial order on $\{1,
\dots, k\}^n$ obtained by the transitive closure of 

$$ (\dots a\dots b\dots)< (\dots b\dots a\dots)\quad \text{ if } 
a>b$$

Then it is easy to show, using (2.4), that $\psi(\e_I)=\e_I+\sum_{I'<I}
c_{I,I'} \e_{I'}$. On the other hand, $I^0(\m)$ is a minimal element
with respect to this order, which proves
$\psi(\e_{I^0(\m)})=\e_{I^0(\m)}$.

To prove that $\psi(PR)_i=(PR)_i^{-1}\psi$, let us first consider the case
$n=2$. Then this identity reduces to  $\Theta\overline{P\bar
C\bar \Theta}=(P\bar C\bar \Theta)^{-1} \Theta$, which is obvious in
view of $\bar C=C^{-1}, \bar \Theta=\Theta^{-1}$. 

For $n\ge 3$, the identity $\psi(PR)_i=(PR)_i^{-1}\psi$ follows by
induction from the $n=2$ case and the fact that $(PR)_i$ is an
intertwiner. Indeed, assume that $W$ is a tensor power of $V$ and
$T:W\to W$ is an intertwining operator such that $\psi
T=T^{-1}\psi$. Then the same is true for the operators $1\o T:V\o W\to
V\o W, T\o 1: W\o V\to W\o V$. This is because the involution $\psi$
on $V\o W$ is given by $\psi=\Theta(\psi\o \psi)$, and therefore
$\psi (1\o T)=\Theta(\psi\o (\psi T))= \Theta(\psi\o
(T^{-1}\psi))= (1\o T^{-1}) \Theta(\psi\o \psi)=(1\o T^{-1})\psi$.
\qed\enddemo

\proclaim{Theorem 2.5} 
Under the assumptions of Proposition~\rom{2.1}, the isomorphism $\phi_\m$
maps the Kazhdan-Lusztig basis in $M$ defined in Theorem~\rom{1.2}
\rom{(}for $u=-1$\rom{)} to the canonical basis in $V^{\o n}[\m]$: 

$$\phi_\m (C_\sigma)=(-1)^{l(\sigma)}b_{\sigma(I^0(\m))}$$
\endproclaim
\demo{Proof} Immediately follows from the definitions and from the
previous proposition.
\qed   \enddemo

So far, we have discussed the relation between the basis $C_\sigma$ in
$M$ defined for $u=-1$ and the canonical basis for $\Uk$. It turns out
that the basis $C_\sigma$ for $u=v^{-2}$ also admits a nice
interpretation in terms of representations of $\Uk$: it is related
with the so-called dual canonical basis.

Let us define on $\Uk$ another structure of Hopf algebra by 

$$\aligned
&\bar\Delta E_i= E_i\o 1+ v^{-H_i}\o E_i,\\
&\bar\Delta F_i= F_i\o v^{H_i}+1 \o F_i,\\
&\bar\Delta v^{H_i}=v^{H_i}\o v^{H_i}
\endaligned 
\tag 2.12
$$

The universal $R$-matrix for this comultiplication is given by 

$$
\bar \Cal R=C\Theta.\tag 2.13
$$

  From now on, let us consider the action of $\Uk$ on tensor powers of
$V$ given by iterations of $\bar \Delta$. Note that since 
$\bar\Delta (v^{H_i})=\Delta (v^{H_i})$, the  weight decomposition for both
actions coincide. Define an antilinear  involution $\psi'$  on tensor
powers of $V$ by the same formulas as (2.10) but replacing $\Theta$ by
$\bar\Theta$, which acts on tensor powers of $V$ via $\bar\Delta$. 

Then the propositions 2.1--2.4 have their analogue. We formulate the
corresponding statements below, using the same notations and
conventions as much as possible. 

\proclaim{Proposition 2.1${}'$} 
In the notations of Proposition~\rom{2.1}, we have 
\roster \item

$$
P\bar R_{V\o V}(\e_a\o \e_b)= v^{1/n}
	\cases \e_b\o \e_a +(v^{-1}-v) \e_a\o \e_b,\quad& a<b,\\
		\e_b\o \e_a, \quad& a>b,\\
		v^{-1} \e_a\o \e_b,\quad& a=b	
	\endcases
\tag 2.14
$$

\item 

The map 
$$
T_i\mapsto v^{-1-1/n} (P\bar R)_i
\tag 2.15
$$
defines an action of the Hecke algebra $\Cal H(S_n)$  on
$V^{\o n}$. This action commutes with the action of $\Uk$; in
particular, it preserves the weight subspaces.

\item 
For $\m\in \Z_+^k, \sum m_i=n$ the map

$$\aligned
\phi'_\m:M&\to V^{\o n}\\
	m_\sigma&\mapsto \e_{\sigma(I^0(\m))}
\endaligned
\tag 2.16$$
is an isomorphism of $\Cal H$ modules $M$ \rom{(}with
$u=v^{-2}$\rom{)} and $V^{\o n}[\m]$.
\endroster
\endproclaim

\proclaim{Proposition 2.3${}'$({\rm \cite{FK}})}
\roster\item
For every $I\in \{1, \dots, k\}^n$ there
exists a unique element $b^I\in V^{\o n}$ such that 

$$\gathered
\psi'(b^I)=b^I,\\
b^I-\e_I\in \bigoplus_{I'} v^{-1}\Z[v^{-1}]\e_{I'},
\endgathered
\tag 2.17
$$
where the sum is taken over all $I'$ having the  same weight as $I$. 

\item Denote by $\<,\>$ the $\Q(v)$-bilinear form on $V^{\o n}$
defined by 

$$\<\e_I, \e_{I'}\>=\delta_{I, w_0(I')},\tag 2.18
$$
where $w_0$ is the longest element in $S_n$: $w_0(1,\dots,
n)=(n,\dots, 1)$. Then:

$$\<b^I, b_{I'}\>=\delta_{I,w_0(I')}.
\tag 2.19
$$
\endroster

The elements $b^I$ form a basis in $V^{\o n}$ which is called the
{\bf dual  canonical basis}. 
\endproclaim

This is the only proposition which requires a separate proof. It does
not immediately follow from Proposition~2.3 since we have replaced $v$
by $v^{-1}$ in all the formulas but left condition (2.11)
unchanged. To prove the proposition, define the elements $b^I$ by
(2.19); then the second condition in (2.17) is satisfied
automatically, and one only needs to check $\psi'(b^I)=b^I$, which is
sufficient to check for $n=2$. We refer the reader to \cite{FK} for
the details  in the case of $\sltwo$; the general case is proved in the
same way.

\proclaim{Proposition 2.4${}'$} The isomorphism $\phi'_\m$ constructed in
Proposition~\rom{2.$1'$} identifies the bar involution in $M$
\rom{(}for $u=v^{-2}$\rom{)} with the involution $\psi'$ in $V^{\o
n}[\m]$.
\endproclaim

\proclaim{Theorem 2.5${}'$} 
Under the assumptions of Proposition~\rom{2.1}, the isomorphism $\phi_\m$
maps the Kazhdan-Lusztig basis in $M$ defined in Theorem~\rom{1.2}
\rom{(}for $u=v^{-2}$\rom{)} to the dual canonial basis in $V^{\o
n}[\m]$:  

$$\phi'_\m (C_\sigma)=b_{\sigma(I^0(\m))}$$
\endproclaim

Therefore, we see that the picture is quite symmetric: we have two
actions of Hecke algebra on the module $M$, which give rise to two
different bases $C_\sigma$. Similarly, we have two different actions
of $\Uk$ on $V^{\o n}$, which give rise to the canonical basis and the
dual canonical basis. 

In particular, let us consider the case $k=n$, and $\m=(1,\dots,
1)$. In this case $J=\varnothing$, and we get the following theorem:

\proclaim{Theorem 2.6} The canonical basis $b_I$ and the dual
canonical basis $b^{I}$ in the zero-weight
subspace of $V^{\o n }$, where $V$ is the fundamental representation of
$U_q(\frak{sl}_n)$, are  given by 

$$
\align
b_{w(n,\dots,1)}=&
   \sum_{y\le w} v^{l(y)-l(w)}  
	\overline{P_{y,w}}\e_{y(n,\dots, 1)} \tag 2.20\\
b^{w(n,\dots,1)}=&
   \sum_{y\le w}(-v)^{l(y)-l(w)}   \overline{P_{y,w}}\e_{y(n,\dots, 1)},
\tag 2.21
\endalign
$$
where $P_{y,w}=P_{y,w}(q), q=v^{-2}$ are the Kazhdan-Lusztig
polynomials for $S_n$.
\endproclaim

Therefore, the problem of calculating Kazhdan-Lusztig polynomials can
be considered as a special case of more general problem of calculating
the canonical basis in a tensor power of the fundamental
representation of $\Uk$. This latter problem, in general, is no
easier than the original one, and there is little hope that it can be
solved in full generality by elementary methods. However, there is one
special case which in which the answer is known: this is the case of
$\U2$, which we will consider in the next section.

\head 3. Canonical basis for $\U2$ and Kazhdan-Lusztig polynomials for
Grassmanians.
\endhead

In this section we recall the known results about the canonical basis
for tensor product of representations of $\U2$ and use the results of
the previous section to transform these results into some statements
about Kazhdan-Lusztig polynomials. It turns out that in this way we
exactly recover the known formulas for Kazhdan-Lusztig polynomials for
Grassmannians, i.e. those corresponding to maximal parabolic subgroups
in $S_n$. These formulas were first obtained by Lascoux and
Sch\"utzenberger \cite{LS} by purely combinatorial methods and later
by Zelevinsky \cite{Z} using a small resolution of singularities of
the corresponding Schubert varieties.

In this section we consider the quantum group $\U2$. It is generated
by the elements $E=E_1, F=F_1, v^H=v^{H_1}$ with the usual commutation
relations and with the comultiplication given by (2.2). The
fundamental representation $V$ is two dimensional: $V=\Q(v)\e_+\oplus
\Q(v)\e_-$, and the action of $\U2$ is given by 

$$\aligned
&E\e_+=0,\quad E\e_-=\e_+,\\
&F \e_+=\e_-,\quad F \e_-=0,\\
&v^H\e_\pm=v^{\pm 1}\e_\pm.
\endaligned
\tag 3.1
$$

The action of the universal $R$-matrix in $V\o V$  is given by

$$\aligned
PR(\e_+\o \e_+ )&=v^{1/2}\e_+\o \e_+,\\
PR(\e_-\o \e_- )&=v^{1/2}\e_-\o \e_-,\\
PR(\e_- \o \e_+ )&=v^{-1/2}\e_+\o \e_-,\\
PR(\e_+ \o \e_- )&=v^{-1/2}\e_-\o \e_++(v^{1/2}-v^{-3/2})\e_+\o \e_-.
\endaligned
\tag 3.2
$$

The relation with the notations of the previous section is given by
$\e_+=\e_1, \e_-=\e_2$.  

In this case it is possible to give an explicit construction of the
dual canonical basis in the tensor power $V^{\o n}$ (and, in fact, in
any tensor product of finite-dimensional representations of
$\U2$). This was done in \cite{FK}. We give here the answer in the
case of interest for us.

Let 
$$a=\e_+\o \e_- -v^{-1}\e_-\o \e_+\in V\o V.
\tag 3.3
$$
One easily check that $a$ is an invariant: $\bar \Delta
E(a)=\bar \Delta(F)(a)=0$ and that $\psi'(a)=a$. Thus, $a$ is an
element of the dual canonical basis in $V^{\o 2}[0]$.  

The following theorem describes the dual canonical basis in tensor
powers $V^{\o m}$. This basis is naturally indexed by sequences $I$ of
pluses and minuses of length $m$ (as before, to relate to the 
notation of the previous section, replace plus by $1$ and minus
by $2$). If $I_1, I_2$ are sequences of length $l,k$ respectively then we
will denote by $I_1|I_2$ the sequence of length $k+l$ obtained by
appending the sequence $I_2$ to the  sequence $I_1$.

\proclaim{Theorem 3.1} The dual canonical basis $b^I$ in $V^{\o m}$ is
given by the following rules:

\roster\item If $I=-|I_1$ then $b^I=\e_-\o b^{I_1}.$

\item If $I=I_1|+$ then $b^I=b^{I_1}\o \e_+.$ 

\item If $I=I_1|+-|I_2$, where $I_1, I_2$ are sequences of length
$i, m-i-2$ respectively then 

$$b^I=a_{i+1 i+2}b^{I_1|I_2},\tag 3.4
$$
where the operator $a_{i+1 i+2}:V^{\o (m-2)}\to V^{\o m}$ is defined
by 

$$a_{i+1 i+2}: \e_{j_1}\o \dots\o \e_{j_{m-2}}\mapsto 
	\e_{j_1}\o\dots\o \e_{j_i}\o a \o \e_{j_{i+1}}\o\dots\o \e_{j_{m-2}}.
\tag 3.5
$$
Here $a\in V\o V$ is given by \rom{(3.3)}.
\endroster\endproclaim

As was mentioned above, this theorem is a special case of the results
in \cite{FK, Section~2.3}; however, this case is not 
difficult to prove from the definitions. 
We leave it as an exercise for the reader. This theorem can be very
neatly visualized usung the graphical calculus; this interpretation
can be found in \cite{FK}. 

\example{Example 3.2} 
$$
\split
b^{++--}=&\e_+\o a\o \e_- - v^{-1}\e_-\o a\o e_+\\ 
	=&\e_{++--}-v^{-1}\e_{-+-+}-v^{-1}\e_{+-+-}+v^{-2}\e_{--++}.
\endsplit
$$
\endexample

Since we have proved that the dual canonical basis in $V^{\o n}$
coincides with the parabolic Kazhdan-Lusztig basis $C_\sigma$ defined
in Theorem~1.2 (for $u=v^{-2}$), the above theorem also gives a very
explicit way to calculate the latter basis for the maximal parabolic
subgroup in $S_n$. In particular, we have the following corollary.

\proclaim{Corollary 3.3} In the notations of Section~\rom{1}, assume that
$J$ is a maximal parabolic subgroup in $W=S_n$. Then each 
parabolic Kazhdan-Lusztig polynomial $P^J_{\tau\sigma}$  \rom{(}for
$u=v^{-2}$\rom{)} is either zero or a power of $q$. 
\endproclaim

\demo{Proof} Immediately follows from Theorems~3.1, 2.5${}'$ and the explicit
formula (3.3) for $a$.
\qed\enddemo

We next describe the canonical basis in tensor powers $V^{\o m}.$ 
In the same way as for the dual canonical basis, we will index the
basis vectors by the length $m$ sequences of pluses and minuses. 
For a sequence $J$ of pluses and minuses denote by $l(J)$ the length
of $J$ and by $J_+$ the number of pluses in $J.$ 

Let $p_n:V^{\o n}\to V^{\o n}$ be the operator of the projection 
onto the unique irreducible $(n+1)$-dimensional submodule of $V^{\o
n}$ (relative to the comultiplication $\Delta$). 
 This operator is sometimes called {\bf the Jones-Wenzl
projector}.  Coordinatewise $p_n$ is just the
$q$-symmetrisation: 
let $I=i_1\dots i_n$ be a sequence which contains  $k$ pluses and 
$n-k$ minuses (in any order). Then 

$$p_n(\e_{i_1}\o \dots\o \e_{i_n})=
\biggl[{
\matrix
n \\    k
\endmatrix}
\biggr]^{-1}
\sum_{J=j_1\dots j_n,J_+=k}
v^{k(n-k)-\lbrace I\rbrace-\lbrace J\rbrace }\e_{j_1}\o \dots\o \e_{j_n}
\tag 3.6
$$

The sum in the right hand side is over all sequences of length $n$ 
that have the same number of pluses as the sequence $i_1\dots i_n,$
and $\lbrace J\rbrace$ is the number of pairs $(a,b), 1\le a<b \le n$
such that $j_a=-$ and $j_b=+$ (compare with Proposition~2.2).  

The quantum binomial coefficient is defined by 
$$
\biggl[{\matrix n \\    k\endmatrix}\biggr]=
\frac{[n]!}{[k]![n-k]!},
$$
where $[n]!=[1][2]\dots[n]$ and $[i]=\frac{v^i-v^{-i}}{v-v^{-1}}.$ 

Also let $p_{i,j}:V^{\o m}\to V^{\o m}$ be given by 
$$p_{i,j}=1^{\o (i-1)}\o p_{j-i+1}\o 1^{\o (m-j)}\tag{3.7}$$ 

\proclaim{Theorem 3.4} The canonical basis $b_I$ in $V^{\o m}$ 
is given by the following rules: 

\roster\item If $I=-|I_1$ then $b_I=\e_-\o b_{I_1}.$

\item If $I=I_1|+$ then $b_I=b_{I_1}\o \e_+.$ 

\item Let  $I=I_1|I_+|I_-|I_2$, where $I_+$ is made entirely of pluses, 
$I_-$ entirely of minuses, and the following conditions are satisfied:

-- $l(I_+)=k>0$, $l(I_-)=l>0$

-- $I_1$ is either empty or ends with  at least $k$ minuses.

-- $I_2$ is either empty or starts with at least $l$ pluses.

Then 
$$
b_I=
\biggl[ {\matrix k+l \\ l \endmatrix}\biggr]
p_{l(I_1)+1,l(I_1)+k+l}
b_{I_1|I_-|I_+|I_2},\tag 3.8
$$

\endroster\endproclaim

This theorem is proved in \cite{K}. The proof is based on Theorem~3.1,
describing the dual canonical basis, and the identity $\<b_I,
b^{I'}\>=\delta_{I, w_0(I')}$ (see (2.19)).

This theorem describes each canonical basis vector as a certain
composition of Jones-Wenzl projectors applied to the vector
$\e_-\o\dots\o \e_-
\o \e_+\o \dots \o \e_+.$ 

\example{Example} $b_{+-++-}=[2][4]p_{1,2}p_{2,5}\e_{--+++}.$ 
\endexample
Note that for some canonical basis vectors such a 
representation is not unique: 

$$b_{+-+-}=[2][3]p_{1,2}p_{2,4}\e_{--++}=[2][3]p_{3,4}p_{1,3}\e_{--++}.$$

By Theorem 2.4, this theorem also describes the Kazhdan-Lusztig basis
in $M$ for $u=-1.$ Thus, in this case the Kazhdan-Lusztig basis
vectors can be expressed as compositions of Jones-Wenzl projectors.

Next we show that 
if the recursive formulas of Theorem~3.4 are written coordinatewise,
one obtains the recursive formulas of 
A.~Lascoux, M.-P.~Sch\"utzenberger \cite{LS} and 
Zelevinski \cite{Z} for 
Kazhdan-Lusztig polynomials in the Grassmannian case.

For a pair of sequences $I,J\in \lbrace +,-\rbrace ^n$ denote by 
$c(I,J)$ the coefficient of $\e_J$ in the decomposition of 
$b_I$ in the product basis of $V^{\o n}.$ 
Order the set $\lbrace +, -\rbrace $ by $+ > -.$ 
Let $I=i_1\dots i_n$ and 
$J=j_1\dots j_n.$ Theorem~3.4 implies that if  
$i_a\ge i_{a+1}$ then 
$p_{a,a+1}b_I =b_I.$ Therefore, if $i_a \ge i_{a+1}, $ 
and $j_a=-, j_{a+1}=+,$ then 
$$
c(I,J)=v^{-1}c(I,J^{(a)}),
\tag 3.9
$$
where $J^{(a)}$ is the sequence obtained by interchanging the $a$th and 
$(a+1)$st elements of $J.$ Thus we can reduce the computation of 
$c(I,J)$ to the case $j_a\ge j_{a+1}$ whenever $i_a\ge i_{a+1}.$ 
If $J$ has this property, we say that $J$ is controlled by $I.$

Before we proceed, we introduce new notations to match the ones in 
\cite{Z}. We can encode $I$ by a sequence of numbers
$b_0,a_1,b_1,a_2,\dots ,b_{m-1},a_m$ as follows. $I$ starts with 
$b_0$ pluses, followed by $a_1$ minuses, followed by $b_1$ pluses, 
etc. $I$ fixed, we encode $J$ (where $J$ is controlled by $I$)
relative to $I$ by  a sequence $x_0,\dots x_{m-1}.$ Namely, $J$ begins 
with $x_0$ pluses, followed by $b_0+a_1-x_0$ minuses, followed by 
$x_1$ pluses, followed by $b_1+a_2-x_1$ minuses, etc. 

We want a recursive formula for $c(I,J)$ in terms of these encodings
of $I$ and $J$. We apply Theorem~3.4. Suppose we can find $i$ such
that $a_i\ge b_i$ and $a_{i+1}\le b_{i+1}$ (this is the most
interesting case; in the notations of Theorem~3.4 it corresponds to
case when both $I_1, I_2$ are non-empty ).  Denote by $L$ the sequence
obtained by interchanging two subsequences of $I$: the subsequence of
$b_i$ pluses and the adjacent subsequence of $a_{i+1}$ minuses. $L$ is
encoded by $b_0, a_1,\dots , b_{i-1},a_{i}+a_{i+1},b_i+b_{i+1},
a_{i+2}, b_{i+2},
\dots , a_m.$

By Theorem 3.4 
$$b_I = 
\biggl[ {\matrix a_{i+1}+b_i \\ b_i \endmatrix}\biggr]
p_{k+1,k+a_{i+1}+b_i}b_L
=\biggl[ {\matrix a_{i+1}+b_i \\ b_i \endmatrix}\biggr]
p_{k+1,k+a_{i+1}+b_i}
\sum_{J'}c(L, J')\e_{J'},\tag 3.10
$$
where $k=b_0+a_1+\dots +b_{i-1}+a_i$. 

Since the projector $p_{k+1,k+1+a_{i+1}+b_i}$ can only change the indices
of $\e_{J'}$ which are in the interval $[k+1,k+1+a_{i+1}+b_i]$, the only
terms in (3.10) which give non-zero contribution to $c(I,J)$ are those
with $J'$ coinciding with $J$ outside of this interval. Using (3.9),
we replace each such $J'$ by a sequence controlled by $L$, which
together with the explicit formula (3.6) for the projector gives 

$$c(I,J)=
\sum_s 
\biggl[ {\matrix x_i  \\ s  \endmatrix} \biggr] 
\biggl[ {\matrix 
 a_{i+1} +b_i-x_i \\ a_{i+1}-s  \endmatrix} \biggr] v^{f(s)}c(L,J_s)
\tag{3.11}$$
where  $J_s$ is the sequence, controlled by $L,$ 
and encoded, relative to $L,$ by 
$$x_0, \dots, x_{i-2}, x_{i-1}+s, x_{i}+x_{i+1}-s, x_{i+2}, \dots
x_{m-1},$$
and
$$f(s)=s(b_i-x_i)+s^2-s(a_i+b_{i-1}-x_{i-1})-x_{i+1}(b_i-x_i+s).
\tag{3.12}$$
Introducing new parameters 
$$\gathered
c_0=x_0, c_{j+1}= c_j+b_j-x_j,\quad 1\le j\le m,  \\
d = c_i-s
\endgathered
$$
we get 
$$c(I,J)=
\sum_s 
\biggl[ {\matrix b_i-c_{i+1}+c_i  \\ c_i-d  \endmatrix} \biggr] 
\biggl[ {\matrix  a_{i+1} +c_{i+1}-c_i \\ c_{i+1}-d  \endmatrix} \biggr]
v^{f(s)}
c(L,J_s)\tag{3.13}$$
and $f(s)$ can be rewritten as 
$$f(s)=(c_i-d)b_i + (c_{i+1}-d)a_{i+1} + g(s)$$
where 
$$\aligned
& g(s)=-((c_i-d)(a_i+b_i)+(c_{i+1}-d)(b_{i+1}+a_{i+1})) + \\
& c_{i-1}c_i-c_i^2 +c_ic_{i+1}-c_{i+1}^2+c_{i+1}c_{i+2}-
  dc_{i-1}-dc_{i+2}+d^2 
\endaligned
$$

Now let 
$$h(I,J)= -\sum_j c_jc_{j+1}+\sum_j c_j(c_j+a_j+b_j)\tag{3.14}$$
and  
$$c^0(I,J)= v^{h(I,J)}c(I,J).\tag{3.15}$$
Then the formula (3.13) becomes 
$$c^0(I,J)=
\sum_s 
\biggl[ {\matrix  b_i-c_{i+1}+c_i  \\ c_i-d  \endmatrix} \biggr] 
\biggl[ {\matrix  a_{i+1} +c_{i+1}-c_i \\ c_{i+1}-d  \endmatrix} \biggr]
v^{(c_i-d)b_i+(c_{i+1}-d)a_{i+1}}
c^0(L,J_s).\tag{3.16}$$

This formula coincides with Zelevinski's inductive formula (\cite{Z}, 
Theorem 2) for the coefficients of Kazhdan-Lusztig polynomials 
for Grassmannians, where his variable $t$ is our $v.$ Note that  
Zelevinski's formula has a power of $v$ different from ours, due to his use 
of shifted quantum binomials, while our quantum binomials are
balanced. 

Observe that formula (3.16) for Kazhdan-Lusztig polynomials
is local.  
Formulas given by Theorems~3.4 and 2.5 for the Kazhdan-Lusztig 
basis have the advantage of 
being global and they make obvious the central role played by the quantum 
group $\U2$ and its intertwiners for the Kazhdan-Lusztig theory in 
the Grassmannian case.  

\head Acknowledgements
\endhead
The works of the authors was partially supported by National Science
Foundations Grants \#DMS-9700765 (first author), \#DMS-9610201 (second
author).

A.K. would like to thank G.~Lusztig for fruitful discussions, and
M.K. is grateful to V.~Protsak for interesting discussions.

\enddocument
\end